\documentclass[twoside,12pt,leqno]{amsart}
\usepackage{amssymb,latexsym,verbatim,pstricks} 

\theoremstyle{plain}
\newtheorem{theorem}{Theorem}

\theoremstyle{definition}

\theoremstyle{definition}              

\renewcommand{\geq}{\geqslant}
\renewcommand{\leq}{\leqslant}

\begin{document}

\newcommand{\C}{\textup{C}}
\newcommand{\CW}{\mathcal{CW}}
\newcommand{\eps}{\varepsilon}
\newcommand{\Farey}{\mathcal{F}}
\newcommand\GL{\textup{GL}}
\newcommand\n{\newline}
\newcommand{\N}{\mathbb{N}}
\newcommand{\p}{\mathfrak{p}}
\newcommand{\Q}{\mathbb{Q}}
\newcommand{\R}{\mathbb{R}}
\newcommand{\SB}{\mathcal{SB}}
\def\u#1{\overbrace{\textstyle{#1}}}
\def\v#1{\vskip#1mm}
\newcommand{\Z}{\mathbb{Z}}

\hyphenation{Frob-enius}

\begin{center}\large\sffamily\mdseries
  Enumerating the rationals from left to right
\end{center}

\title[\tiny\upshape\rmfamily Enumerating the rationals from left to right]{}


\author{\sffamily S.\,P. Glasby}

\address[]{\kern -4.4mm Department of Mathematics\\   
Central Washington University\\
WA 98926-7424, USA.\newline {\tt http:/\kern-1pt/www.cwu.edu/$\sim$glasbys/}}


\maketitle

\vskip-2mm
\centerline{\noindent\Small 2010 Mathematics subject classification: 11B57, 11B83, 11B75}


There are three well-know sequences used to enumerate the rationals: the
Stern-Brocot sequences $\SB_n$, the Calkin-Wilf sequences $\CW_n$,
and the Farey sequences $\Farey_n$. The purpose of this note is to show that
all three sequences can be constructed (left-to-right) using almost identical
recurrence relations. The Stern-Brocot (S-B) and
Calkin-Wilf (C-W) sequences give rise to complete binary trees related
to the following rules:
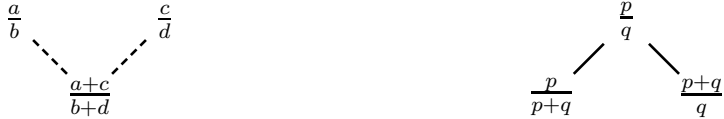
\v2
\begin{center}
  \begin{pspicture}(0,0)(2,1)
  \rput(0,1){$\frac ab$}
  \rput(2,1){$\frac cd$}
  \rput(1,0){$\frac{a+c}{b+d}$}
  \psline[linewidth=1pt,linestyle=dashed,dash=3pt 2pt](0.25,0.75)(0.7,0.3)
  \psline[linewidth=1pt,linestyle=dashed,dash=3pt 2pt](1.3,0.3)(1.75,0.75)
  \end{pspicture}
  \hskip50mm
  \begin{pspicture}(0,0)(2,1)
  \rput(1,1){$\frac pq$}
  \rput(0,0){$\frac p{p+q}$}
  \rput(2,0){$\frac {p+q}q$}
  \psline[linewidth=1pt](0.3,0.3)(0.7,0.7)
  \psline[linewidth=1pt](1.3,0.7)(1.7,0.3)
  \end{pspicture}
\end{center}
\v1
These trees have many beautiful algebraic, combinatorial,
computational, and geometric properties~\cite{CW,M,H06}.
Well-written introductions to the S-B tree and Farey sequences
can be found in~\cite{GKP}, and to the C-W tree in~\cite{CW}.
We shall focus on {\it sequences} rather than {\it trees}.

Two fractions $\frac ab<\frac cd$ are called {\it adjacent} if $bc-ad=1$.
Adjacent fractions are necessarily reduced, i.e. $\gcd(a,b)=\gcd(c,d)=1$.
The {\it mediant} of $\frac ab<\frac cd$ is $\frac{a+c}{b+d}$.
Simple algebra shows if $\frac ab<\frac cd$ are adjacent, then
$\frac ab<\frac{a+c}{b+d}<\frac cd$ are pairwise adjacent (and hence reduced).
The sequences $\SB_n$ are defined recursively: $\SB_0=[\frac01,\frac10]$
represents $0$ and $\infty$ as fractions, and $\SB_n$ is computed from
$\SB_{n-1}$ by inserting mediants between consecutive fractions. Thus
\[
  \SB_1=\left[\frac01,{\bf\frac11},\frac10\right],\quad
  \SB_2=\left[\frac01,{\bf\frac12},\frac11,{\bf\frac21},\frac10\right],\quad
  \SB_3=\left[\frac01,{\bf\frac13},\frac12,{\bf\frac23},\frac11,
    {\bf\frac32},\frac21,{\bf\frac31},\frac10\right],\ \dots.
\]
A simple induction shows that $|\SB_n|=2^n+1$. Thus $2^{n-1}$ mediants
are inserted into $\SB_{n-1}$ to form $\SB_n$. The C-W sequences are defined
using the right rule above:
\[
  \CW_1:=\left[\frac11\right],\ 
  \CW_2=\left[\frac12,\frac21\right],\ 
  \CW_3=\left[\frac13,\frac32,\frac23,\frac31\right],\ 
  \CW_4=\left[\frac14,\frac43,\frac35,\frac52,\frac25,
    \frac53,\frac34,\frac41\right], \dots.
\]
A simple induction shows that $|\CW_n|=2^{n-1}$. Another simple induction
(see~\cite[p360]{CW}) shows that the fractions in $\CW_n$ have the form
\[
  \CW_n=\left[\frac{b_{-1}}{b_0},\frac{b_0}{b_1},\dots,
    \frac{b_{N-2}}{b_{N-1}}\right]\qquad\textup{where $N=2^{n-1}$,}
\]
and the denominator of a given fraction is the numerator
of the succeeding fraction. Indeed, this property obtains even when the
sequences $\CW_1,CW_2,CW_3,\dots$ are concatenated to form
$\CW_\infty:=
  [\frac1{\bf1},\u{\frac{\bf1}2,\frac2{\bf1}},
  \u{\frac{\bf1}3,\frac32,\frac23,\frac3{\bf1}},
  \u{\frac{\bf1}4,\dots,\frac4{\bf1}},\dots]
  =\left[\frac{c_0}{c_1},\frac{c_1}{c_2},\frac{c_2}{c_3},\dots\right]$.

The Farey sequence of order~$n$ contains all the reduced fractions
$\frac pq$ with $0\leq p\leq q\leq n$, in their natural order. Thus
\[
  \Farey_1=\left[\frac01,\frac11\right],\ 
  \Farey_2=\left[\frac01,{\bf\frac12},\frac11\right],\ 
  \Farey_3=\left[\frac01,{\bf\frac13},\frac12,{\bf\frac23},\frac11\right],\ 
  \Farey_4=\left[\frac01,{\bf\frac14},\frac13,\frac12,\frac23,{\bf\frac34},
    \frac11\right],\dots. 
\]
A standard way to compute $\Farey_n$ from
$\Farey_{n-1}$ is to insert mediants between consecutive fractions of
$\Farey_{n-1}$ only when it gives a denominator of size~$n$
(see~\cite[p118]{GKP}). Thus $\Farey_n$ is a subsequence of $\SB_n$.
It is easy to prove that $|\Farey_n|=1+\sum_{j=1}^n\varphi(j)$
where $\varphi(n)$
denotes the number of reduced fractions $\frac an$ with $1\leq a<n$.
The mediant rule above implies that consecutive fractions in
$\SB_n$ and $\Farey_n$ are adjacent (see also~\cite[p119]{GKP}).

It is shown in~\cite{GKP} and~\cite{CW} that $\SB_\infty$ and
$\CW_\infty$ contain {\it every} (reduced) positive rational
{\it precisely once}.
Although $\SB_n, \CW_n, \Farey_n$ are defined ``top-down'' they can be computed
from ``left to right'' via almost identical recurrence relations.

\begin{theorem}\label{T:SB}
Write $\SB_n=\left[\frac{a_{-1}}{b_{-1}},\frac{a_0}{b_0},\frac{a_1}{b_1},
    \dots,\frac{a_{N-1}}{b_{N-1}}\right]$ where $N=2^n$. Then
\begin{subequations}\label{E:SB}
\begin{align}
  &a_{-1}=0,\ a_0=1,&&a_{i}=k_ia_{i-1}-a_{i-2}&&\textup{for $1\leq i<N$},
  \label{E:SBRR1}\\
  &b_{-1}=1,\ b_0=n,&&b_{i}=k_ib_{i-1}-b_{i-2}&&\textup{for $1\leq i<N$},
  \label{E:SBRR2}
\end{align}
\end{subequations}
where $k_i=2\log_2 |i|_2 \,+\,1$, and $|i|_2$ denotes the largest power of~$2$
dividing~$i$.
\end{theorem}

\begin{theorem}\label{T:CWRR}
Write $\CW_\infty=
  \left[\frac{a_0}{a_1},\frac{a_1}{a_2},\dots,\frac{a_{i-1}}{a_i},
  \dots\right]$
and $\CW_n=\left[\frac{b_{-1}}{b_0},\frac{b_0}{b_1},\frac{b_1}{b_2},
    \dots,\frac{b_{N-2}}{b_{N-1}}\right]$ where $N=2^{n-1}$.
Then the $a_i$ and $b_i$ can be
computed via the recurrence relations
\begin{subequations}\label{E:CWRR}
\begin{align}
  &a_{-1}=0,\ a_0=1,&&a_i=k_ia_{i-1}-a_{i-2}&&\textup{for $1\leq i<\infty$},
  \label{E:CWRR2}\\
  &b_{-1}=1,\ b_0=n,&&b_{i}=k_ib_{i-1}-b_{i-2}&&\textup{for $1\leq i<N$},
  \label{E:CWRR1}
\end{align}
\end{subequations}
where $k_i=2\log_2 |i|_2 \,+\,1$.
\textup{[Note that $\nu_2(i):=\log_2 |i|_2$ is the largest $\nu\in\Z$
satisfying $2^\nu\mid i$.]}
\end{theorem}

\begin{theorem}\label{T:FareyRR}
Write the Farey sequence $\Farey_n$ of order~$n$ as
$\Farey_n=[\frac{A_{-1}}{B_{-1}},\frac{A_0}{B_0},\frac{A_1}{B_1},\dots]$.
Then the numerators $A_i$, and the denominators $B_i$ can be computed
via the recurrence relations
\begin{subequations}\label{E:FareyRR}
\begin{align}
  &A_{-1}=0,\ A_0=1,&&A_{i}=K_iA_{i-1}-A_{i-2}
    &&\textup{for $1\leq i<N$},\label{E:FRR1}\\
  &B_{-1}=1,\ B_0=n,&&B_{i}=K_iB_{i-1}-B_{i-2}
    &&\textup{for $1\leq i<N$},\label{E:FRR2}
\end{align}
\end{subequations}
where $K_i=\left\lfloor\frac{B_{i-2}+n}{B_{i-1}}\right\rfloor$,
and $N=\sum_{j=1}^n\varphi(j)$.
\end{theorem}

To illustrate Theorem~\ref{T:SB},
$\SB_4$ can be computed from left to right using the table
\begin{center}
  \begin{tabular}{|c|c|c|c|c|c|c|c|c|c|c|c|c|c|c|c|c|c|} \hline
  $i$&$-1$&0&1&2&3&4&5&6&7&8&9&10&11&12&13&14&15\\ \hline
  $a_i$&0&1&1&2&1&3&2&3&1&4&3&5&2&5&3&4&1\\ \hline
  $b_i$&1&4&3&5&2&5&3&4&1&3&2&3&1&2&1&1&0\\ \hline
  $k_i$& & &1&3&1&5&1&3&1&7&1&3&1&5&1&3&1\\ \hline
  \end{tabular}\,.
\end{center}
The numbers $k_i$ are the same as the numbers $k'_i$ generated by the recurrence
$k'_1=1$, $k'_{2j+1}=1$, $k'_{2j}=2k'_j+1$ for $j\geq0$. (Proof by induction:
$k_1=k'_1$ and $k_{2j+1}=1$, $k_{2j}=2k_j+1$ hold for $j\geq1$
as $|2j+1|_2=1$ and $|2j|_2=2|j|_2$.
Thus $k_i=k'_i$ for all $i\geq1$.)

\begin{proof}[{\bf Proof} \textup{(of Theorem~\ref{T:SB}).}]
Our proof uses induction on $n$. It suffices to prove (\ref{E:SBRR1}) as the
proof of (\ref{E:SBRR2}) is similar (just change the $a$'s to $b$'s). Clearly
(\ref{E:SBRR1}) is true for $n=0$ as $\SB_0=[\frac01,\frac10]$.
Assume $n>0$ and (\ref{E:SBRR1}) is true for $\SB_{n-1}$.
Let $\frac{a'_{-1}}{b'_{-1}},\frac{a'_{0}}{b'_{0}},\dots,
\frac{a'_{N/2-1}}{b'_{N/2-1}}$ be the fractions in $\SB_{n-1}$.
The way mediants are inserted to create $\SB_n$ is shown below:
\v1
\begin{center}
  \begin{pspicture}(10,2.3)
  \rput(0,2){$\SB_{n-1}$}\rput(0,0){$\SB_n$}
  \rput(2,2){$\frac{a'_{j-2}}{b'_{j-2}}$}
  \rput(6,2){$\frac{a'_{j-1}}{b'_{j-1}}$}
  \rput(10,2){$\frac{a'_j}{b'_j}$}
  \rput(2,0){$\frac{a_{2j-3}}{b_{2j-3}}$}
  \rput(4,0){$\bf\frac{\boldsymbol{a}_{2j-2}}{\boldsymbol{b}_{2j-2}}$}
  \rput(6,0){$\frac{a_{2j-1}}{b_{2j-1}}$}
  \rput(8,0){$\bf\frac{\boldsymbol{a}_{2j}}{\boldsymbol{b}_{2j}}$}
  \rput(10,0){$\frac{a_{2j+1}}{b_{2j+1}}$}
  \psline[linewidth=1pt,linestyle=dashed](2.4,1.6)(3.6,0.4)
  \psline[linewidth=1pt,linestyle=dashed](4.4,0.4)(5.6,1.6)
  \psline[linewidth=1pt,linestyle=dashed](6.4,1.6)(7.6,0.4)
  \psline[linewidth=1pt,linestyle=dashed](8.4,0.4)(9.6,1.6)
  \psline[linewidth=1.5pt,linestyle=dotted](2,0.5)(2,1.5)
  \psline[linewidth=1.5pt,linestyle=dotted](6,0.5)(6,1.5)
  \psline[linewidth=1.5pt,linestyle=dotted](10,0.5)(10,1.5)
  \end{pspicture}
\end{center}
\v2
where dotted lines denote the repetition of a fraction, and
dashed lines denote the formation of a mediant. The repetition of
fractions means
\begin{equation}\label{E:repeat}
  a_{2j-1}=a'_{j-1}\quad\textup{and}\quad b_{2j-1}=b'_{j-1}
    \quad\textup{for $0\leq j<N/2$,} 
\end{equation}
and the formation of mediants means
\begin{equation}\label{E:mediant}
  a_{2j}=a'_{j-1}+a'_j\quad\textup{and}\quad b_{2j}=b'_{j-1}+b'_j
    \quad\textup{for $0\leq j<N/2$.}
\end{equation}
We prove (\ref{E:SBRR1}) using induction on $i$. Certainly
(\ref{E:SBRR1}) is true for $i=-1,0$ as $\SB_n$ starts with
$\frac01,\frac1n$. Suppose now that $i\geq1$, and consider the case when
$i$ is even and odd separately.\n
\textsc{Case 1.} $i=2j$ is even and $j\geq1$. The following shows that
(\ref{E:SBRR1}) holds for even $i$:
\begin{align*}
  k_{2j}a_{2j-1}-a_{2j-2}
    &=(k_j+2)a_{2j-1}-a_{2j-2}&&\textup{as $k_{2j}=k_j+2$,}\\
    &=(k_j+2)a'_{j-1}-(a'_{j-2}+a'_{j-1})
      &&\textup{by (\ref{E:repeat}) and (\ref{E:mediant}),}\\
    &=k_ja'_{j-1}-a'_{j-2}+a'_{j-1}&&\textup{canceling $a'_{j-1}$,}\\
    &=a'_j+a'_{j-1}&&\textup{as $a'_j=k_ja'_{j-1}-a'_{j-2}$ by induction,}\\
    &=a_{2j}&&\textup{by (\ref{E:mediant}).}
\end{align*}
\textsc{Case 2.} $i=2j+1$ is odd and $j\geq1$. 
Note that $k_{2j+1}=1$, $a_{2j-1}=a'_{j-1}$, and $a_{2j+1}=a'_j$
by (\ref{E:repeat}). These equations and (\ref{E:mediant}) now imply 
\[
  k_{2j+1}a_{2j}-a_{2j-1}=a_{2j}-a_{2j-1}=(a'_{j-1}+a'_j)-a'_{j-1}=a'_j=a_{2j+1}
\]
as desired. This completes the inductions on $i$ and $n$.
\end{proof}

A different (and very interesting) method for computing terms
of $\SB_n$ is given in~\cite{BBT}. It uses continued
fraction expansions and ``normal additive factorizations''.
As the recurrence (\ref{E:SBRR1}) is independent of~$n$, the numerators for
$\SB_{n-1}$ reappear as the first $2^{n-1}+1$ numerators for $\SB_n$.
We now show that (half of) the denominators $b_i$ in $\SB_n$ reappear
(remarkably!) for $\CW_n$, and the numerators $a_i$ also reappear
in $\CW_\infty$. Accordingly, we shall use the {\it same notation}
$a_i,b_i$ in Theorem~\ref{T:CWRR} as in Theorem~\ref{T:SB}.

\begin{proof}[{\bf Proof} \textup{(of Theorem~\ref{T:CWRR}).}]
The following diagram of the C-W tree (with rules)
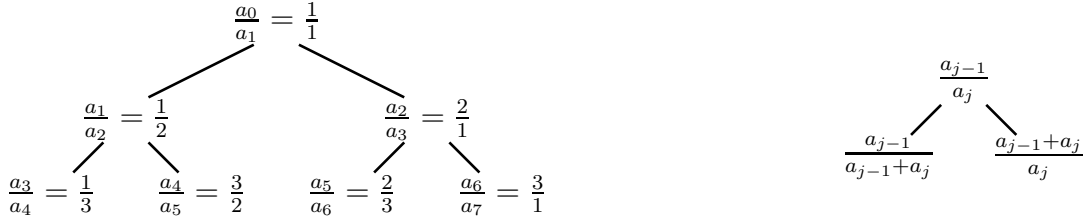
\v2
\begin{center}
  \begin{pspicture}(0,0)(6,2.3)
  \rput(0,0){$\frac{a_3}{a_4}=\frac13$}
  \rput(2,0){$\frac{a_4}{a_5}=\frac32$}
  \rput(4,0){$\frac{a_5}{a_6}=\frac23$}
  \rput(6,0){$\frac{a_6}{a_7}=\frac31$}
  \rput(1,1){$\frac{a_1}{a_2}=\frac12$}
  \rput(5,1){$\frac{a_2}{a_3}=\frac21$}
  \rput(3,2.3){$\frac{a_0}{a_1}=\frac11$}
  \psline[linewidth=1pt](0.3,0.3)(0.7,0.7)
  \psline[linewidth=1pt](1.3,0.7)(1.7,0.3)
  \psline[linewidth=1pt](4.3,0.3)(4.7,0.7)
  \psline[linewidth=1pt](5.3,0.7)(5.7,0.3)
  \psline[linewidth=1pt](1.3,1.3)(2.7,2)
  \psline[linewidth=1pt](3.3,2)(4.7,1.3)
  \end{pspicture}
  \hskip50mm
  \begin{pspicture}(0,0)(2,1.5)
  \rput(1,1.5){$\frac{a_{j-1}}{a_j}$}
  \rput(0,0.5){$\frac{a_{j-1}}{a_{j-1}+a_j}$}
  \rput(2,0.5){$\frac{a_{j-1}+a_j}{a_j}$}
  \psline[linewidth=1pt](0.3,0.8)(0.7,1.2)
  \psline[linewidth=1pt](1.3,1.2)(1.7,0.8)
  \end{pspicture}
\end{center}
\v1
shows that the numbers $a_i$ must satisfy the recurrence relation:
\begin{equation}\label{E:a}
  a_0=1, \quad a_{2j-1}=a_{j-1}\quad\textup{and}\quad
  a_{2j}=a_{j-1}+a_j\quad\textup{for $j>0$.}
\end{equation}
The different recurrence relations (\ref{E:a}) and (\ref{E:CWRR2})
determine the values $a_0,a_1,a_2,\dots$.
We must prove, therefore, that both recurrence relations generate the {\it same}
numbers. For clarity, we write the numbers produced by (\ref{E:a}) as~$a'_i$.
Thus
\begin{equation}\label{E:a'}
  a'_{-1}=0,\quad a'_0=1, \quad
    a'_{2j-1}\overset{(\ref{E:a'}.1)}=a'_{j-1}\quad\textup{and}\quad
  a'_{2j}\overset{(\ref{E:a'}.2)}{=}a'_{j-1}+a'_j\quad\textup{for $j\geq0$.}
\end{equation}
(Note that the definition $a'_{-1}:=0$ is consistent with
$a'_{2j-1}=a'_{j-1}$ and $a'_{2j}=a'_{j-1}+a'_j$ when $j=0$.)
Our goal is to prove $a_i$ defined by (\ref{E:CWRR2}) equals $a'_i$ defined
by (\ref{E:a'}) for $i\geq-1$.

We use induction on $i$. Certainly $a_i=a'_i$ holds for $i=-1,0$.
Assume $i\geq1$ and $a_0=a'_0$, $a_1=a'_1$,\dots, $a_{i-1}=a'_{i-1}$.
Consider the cases when $i$ is odd and even separately.\n
\textsc{Case 1.} $i=2j-1$ where $j\geq1$. Then
\begin{align*}
  a_{2j-1}&=k_{2j-1}a_{2j-2}-a_{2j-3}&&\quad\textup{by (\ref{E:CWRR2}),}\\
         &=a_{2j-2}-a_{2j-3}&&\quad\textup{as $k_{2j-1}=1$,}\\
         &=a'_{2j-2}-a'_{2j-3}&&\quad\textup{by induction on $i$,}\\
         &=a'_{j-2}+a'_{j-1}-a'_{j-2}&&
           \quad\textup{by (\ref{E:a'}.2) and (\ref{E:a'}.1),}\\
         &=a'_{2j-1}&&\quad\textup{by (\ref{E:a'}.1).}
\end{align*}
\textsc{Case 2.} $i=2j$ where $j\geq1$. Then
\begin{align*}
  a_{2j}&=k_{2j}a'_{2j-1}-a_{2j-2}&&\quad\textup{by (\ref{E:CWRR2}) and Case 1,}\\
         &=(k_j+2)a'_{2j-1}-a'_{2j-2}&&
           \quad\textup{by $k_{2j}=k_j+2$ and induction,}\\
         &=(k_j+2)a'_{j-1}-(a'_{j-2}+a'_{j-1})&&
           \quad\textup{by (\ref{E:a'}.1) and (\ref{E:a'}.2),}\\
         &=k_ja'_{j-1}-a'_{j-2}+a'_{j-1}&&\quad\textup{canceling $a'_{j-1}$,}\\
         &=a'_j+a'_{j-1}&&
           \quad\textup{by induction on $i$ and (\ref{E:CWRR2}),}\\
         &=a'_{2j}&&\quad\textup{by (\ref{E:a'}.1).}
\end{align*}
This completes the inductive proof of (\ref{E:CWRR2}).

The proof of (\ref{E:CWRR1}) is now straightforward. As
$\CW_n$ is a subsequence of $\CW_\infty$, there exists an $m$ for which
$\frac{a_{m-1}}{a_m}$ equals the first fraction
$\frac{b_{-1}}{b_0}=\frac1n$ of $\CW_n$. Thus $a_{m-1}=b_{-1}=1$ and
$a_m=b_0=n$.
Since the recurrences (\ref{E:SBRR1}) and (\ref{E:SBRR2}) have the
same form, it follows that $a_{m+i}=b_i$ for $1\leq i<N$. Thus (\ref{E:SBRR2})
holds and (\ref{E:CWRR1}), which is the same, also holds.
\end{proof}

Theorem~\ref{T:FareyRR} is previously known (see Exercise~4.61
in~\cite[p150]{GKP}). We include Theorem~\ref{T:FareyRR} and its proof both
for comparison with Theorems~\ref{T:SB} and~\ref{T:CWRR}, and for the
reader's convenience.

\begin{proof}[{\bf Proof} \textup{(of Theorem~\ref{T:FareyRR}).}]
Our proof uses induction on $i$. As the first two fractions of $\Farey_n$
are $\frac01$ and $\frac1n$, the recurrences (\ref{E:FRR1},b)
are correct for $i=-1,0$. Suppose now that $i>0$ and that (\ref{E:FRR1},b)
are correct for subscripts less than $i$.
Thus $\frac{A_{i-2}}{B_{i-2}}$ and $\frac{A_{i-1}}{B_{i-1}}$ are consecutive
fractions of $\Farey_n$, and we wish to show that the next fraction
is $\frac{A_i}{B_i}$ where $A_i=K_iA_{i-1}-A_{i-2}$ and $B_i=K_iB_{i-1}-B_{i-2}$.
As consecutive Farey fraction are adjacent (i.e they satisfy $bc-ad=1$),
we know by 
induction that $A_{i-1}B_{i-2}-B_{i-1}A_{i-2}=1$. However, the recurrences
(\ref{E:FRR1},b) extend this property as
\begin{align}\label{E:1}
  A_iB_{i-1}-B_iA_{i-1}
    &=(K_iA_{i-1}-A_{i-2})B_{i-1}-(K_iB_{i-1}-B_{i-2})A_{i-1}\\
    &=A_{i-1}B_{i-2}-B_{i-1}A_{i-2}=1.\notag
\end{align}

Consider the inequalities
$\frac{B_{i-2}+n}{B_{i-1}}-1<K_i\leq \frac{B_{i-2}+n}{B_{i-1}}$.
Multiplying by $B_{i-1}$ and subtracting $B_{i-2}$ gives
$n-B_{i-1}<B_i\leq n$. It follows from (\ref{E:1}) and $0<B_i\leq n$ that
$\frac{A_{i-2}}{B_{i-2}}<\frac{A_{i-1}}{B_{i-1}}<\frac{A_i}{B_i}$.
Suppose that $\frac ab$ is the next fraction in $\Farey_n$ after
$\frac{A_{i-1}}{B_{i-1}}$.
Then we know $\frac{A_{i-1}}{B_{i-1}}<\frac ab\leq\frac{A_i}{B_i}$, and
we must show $\frac ab=\frac{A_i}{B_i}$. If not, then 
\begin{equation}\label{E:2}
  A_ib-aB_i\overset{(\ref{E:2}.1)}{\geq}1\qquad\text{and}
  \qquad aB_{i-1}-bA_{i-1}\overset{(\ref{E:2}.2)}{\geq} 1.
\end{equation}
Multiplying (\ref{E:2}.1) by $B_{i-1}$, and (\ref{E:2}.2) by $B_i$, and then
adding gives
\[
  n<B_{i-1}+B_i\leq (A_ib-aB_i)B_{i-1}+B_i(aB_{i-1}-bA_{i-1})
  =(A_iB_{i-1}-B_iA_{i-1})b=b.
\]
This is a contradiction since $\frac ab\in\Farey_n$ has $b\leq n$.
Hence $\frac ab=\frac{A_i}{B_i}$. As both fractions are reduced
(and $a,b,B_i>0$), we conclude that $a=A_i$ and $b=B_i$, as desired.
\end{proof}

The On-Line Encyclopedia of Integer Sequences~\cite{OEIS} has a wealth of
useful information about the sequences $a_0,a_1,a_2,\dots$ (A002487),
and $k_1,k_2,k_3,\dots$ (A037227), however, the connection in
Theorem~\ref{T:CWRR} between these sequences is new.
Note that $a_n$ counts the number of ways that~$n$ can be written as a sum of
powers of~2, each power being used at most twice. For example, $a_4=3$ as
$2^2=2+2=2+1+1$. Finally, we remark that each positive fraction
$\frac pq$ can be associated with a string of $L$'s and $R$'s
denoting its position in a binary tree~\cite[p119]{GKP}.
A simple induction (which we omit) shows
that the S-B string of $\frac pq$ equals the {\it reverse} of the
C-W string of $\frac pq$.

\end{document}